\documentclass[12pt]{article}
\usepackage{amssymb, amsfonts, amsmath, amsthm, url, graphicx}

\def\3{\subset }
\def\4{\subseteq }

\def\0{\leqno}

\def\barr{\begin{array}}
\def\earr{\end{array}}
\def\dd{\displaystyle}

\def\Z{{\rlap{$\kern2pt{\rm Z}$}{\rm Z}\,}}
\def\bld#1#2{{\buildrel{#1}\over{#2}}}
\def\st#1#2{{\mathrel{\mathop{#2}\limits_{#1}}{}\!}}
\def\stb#1#2#3{{\st{{#1}}{\bld{{#2}}{#3}}{}\!}}
\def\xmare#1#2{\stb{#1}{#2}{\mbox{\Huge$\times$}}}

\def\n{\noindent }
\def\bk{\bigskip}
\title{\bf The posets of classes of isomorphic subgroups of finite groups}
\author{Marius T\u arn\u auceanu}
\date{February 17, 2015}

\begin{document}

\maketitle

\begin{abstract}
In this paper we introduce and study the poset of equivalence
classes of subgroups of a finite group $G$, induced by the
isomorphism relation. This contains the well-known lattice of
solitary subgroups of $G$. We prove that in several particular
cases it determines the structure of $G$.
\end{abstract}


poset/lattice isomorphisms.

\section{Introduction}

The relation between the structure of a group and the structure of
its lattice of subgroups constitutes an important domain of
research in group theory. The topic has enjoyed a rapid
development starting with the first half of the 20th century. Many
classes of groups determined by different properties of partially
ordered subsets of their subgroups (especially lattices of
subgroups) have been identified. We refer to Suzuki's book
\cite{12}, Schmidt's book \cite{11} or the more recent book
\cite{14} by the author for more information about this theory.

It is an usual technique to consider an equivalence relation
$\sim$ on an algebraic structure and then to study the factor set
with respect to $\sim$, partially ordered by certain ordering
relations. In the case of subgroup lattices, one of the most
significant example is the poset $C(G)$ of conjugacy classes of
subgroups of a group $G$ (see \cite{2,3,4} and \cite{9,10}). The
current paper deals with the more general equivalence relation on
the subgroup lattice of $G$ induced by isomorphism. It leads
to the set Iso($G$) consisting of all equivalence classes of
isomorphic subgroups of $G$, that becomes a poset under a
su\-i\-ta\-ble ordering relation. Its detailed study is the main
goal of our paper. We investigate the finite groups $G$ for which
the corresponding poset Iso($G$) is a lattice and, in particular,
a chain. We also give some information about the finite groups
$G_1$ and $G_2$ for which the posets Iso($G_1$) and Iso($G_2$) are
isomorphic.

In the following for a finite group $G$ we will denote by $L(G)$
the subgroup lattice of $G$. Recall that $L(G)$ is a complete
bounded lattice with respect to set inclusion, having initial
element the trivial subgroup 1 and final element $G$, and its
binary operations $\wedge, \vee$ are defined by
$$H\wedge K=H\cap K,\ H\vee K=\langle H\cup K\rangle, \mbox{ for all } H,K\in
L(G).$$Two groups $G_1$ and $G_2$ will be called
\textit{L-isomorphic} if their subgroup lattices $L(G_1)$ and
$L(G_2)$ are isomorphic. We also recall that an important modular
sublattice of $L(G)$ is constituted by the normal subgroup lattice
of $G$, usually denoted by $N(G)$.

The paper is organized as follows. In Section 2 we present some
basic properties and results on the poset Iso($G$) associated to a
finite group $G$. A complete description of this poset is given
for several remarkable groups. Section 3 deals with the finite
groups having the same poset of isomorphic subgroups. In the final
section some conclusions and further research directions are
indicated.

Most of our notation is standard and will usually not be repeated
here. Basic definitions and results on lattices and groups can be
found in \cite{1,5} and \cite{6,7,13}, respectively. For subgroup
lattice concepts we refer the reader to \cite{11,12,14}.

\section{The poset Iso($G$)}

Let $G$ be a finite group and ${\rm Iso}(G)$ be the set of
equivalence classes of subgroups of $G$ with respect to the
isomorphism relation, that is
$${\rm Iso}(G)=\{[H] \mid H\in L(G)\}, \mbox{ where }
[H]=\{K\in L(G) \mid K\cong H\}.$$Then it is easy to see that
${\rm Iso}(G)$ can be partially ordered by defining
$$[H_1]\leq [H_2] \mbox{ if and only if } K_1\subseteq K_2 \mbox{ for some } K_1\in [H_1] \mbox{ and } K_2\in
[H_2].$$

We remark that $\leq\hspace{0,5mm}$ is weaker than the usual
ordering relation on $C(G)$ and that the isomorphism relation is
not a congruence on $L(G)$, even if in many cases the poset $({\rm
Iso}(G),\leq)$ becomes a lattice. We also must mention that the
subposet of ${\rm Iso}(G)$ determined by all classes with a unique
element is in fact the lattice ${\rm Sol}(G)$ of solitary
subgroups of $G$, introduced and studied in \cite{8}. \bk

First of all, we will look at the poset ${\rm Iso}(G)$ associated
to some finite groups of small orders.

\bk\n{\bf Examples.}
\begin{itemize}
\item[\rm 1.] ${\rm Iso}(\mathbb{Z}_p)$ is a chain of length 1, for any prime $p$.
\item[\rm 2.] ${\rm Iso}(\mathbb{Z}_p\times\mathbb{Z}_p)\cong{\rm Iso}(\mathbb{Z}_{p^2})$ is a chain of length 2, for any prime $p$.
\item[\rm 3.] ${\rm Iso}(\mathbb{Z}_6)\cong{\rm Iso}(S_3)\cong{\rm Iso}(D_{10})$ is a direct product of two chains of length 2.
\item[\rm 4.] ${\rm Iso}(\mathbb{Z}_2^3)\cong{\rm Iso}(\mathbb{Z}_8)\cong{\rm Iso}(Q_8)$ is a chain of length 3.
\item[\rm 5.] ${\rm Iso}(\mathbb{Z}_2\times\mathbb{Z}_4)\cong{\rm Iso}(D_8)$ is the lattice ${\rm C}_5$ (see page 5 of \cite{11}).
\item[\rm 6.] ${\rm Iso}(A_4)$ is the pentagon lattice ${\rm E}_5$ (see page 5 of \cite{11}).
\end{itemize}
\smallskip

It is well-known that a finite nilpotent group $G$ can be written
as the direct product of its Sylow subgroups
$$G=\xmare{i=1}{k}G_i,$$
where $|G_i|=p_i^{\alpha_i}$, for all $i=1,2,...,k$. Since the
subgroups of a direct pro\-duct of groups having coprime orders
are also direct products (see Corollary of (4.19), \cite{13}, I),
one obtains that
$$L(G)\cong\xmare{i=1}{k}L(G_i).\0(*)$$This lattice direct decomposition
is often used in order to reduce many pro\-blems on $L(G)$ to the
subgroup lattices of finite $p$-groups. We easily observe that
$(*)$ leads to
$${\rm Iso}(G)\cong\xmare{i=1}{k}{\rm Iso}(G_i),\0(**)$$that is the
decomposability of $L(G)$ implies the decomposability of ${\rm
Iso}(G)$. Moreover, all posets ${\rm Iso}(G_i)$, $i=1,2,...,k$,
are indecomposable, since each group $G_i$ possesses a unique
class of isomorphism of subgroups of order $p_i$. The above
example 3 shows that the converse implication fails: ${\rm
Iso}(S_3)$ is decomposable, in contrast with $L(S_3)$.
\bigskip

In the following we shall focus on describing the poset ${\rm
Iso}(G)$ for several important classes of finite groups $G$. We
start with abelian groups, for which we already know that the
study is reduced to abelian $p$-groups.

\bk\n{\bf Proposition 2.1.} {\it Let
$G=\xmare{i=1}{k}\mathbb{Z}_{p^{\alpha_i}}$ be a finite abelian
$p$-group. Then $\hspace{0,5mm}{\rm Iso}(G)$ is an indecomposable
distributive lattice and
$$|{\rm Iso}(G)|\leq\dd\sum_{i=0}^{\alpha}\pi(i),$$where
$\alpha=\alpha_1+\alpha_2+\cdots+\alpha_k$ and $\pi(i)$ denotes
the number of partitions of $\,i$, for all $i=0,1,...,\alpha$.}
\bigskip

\n{\bf Proof.} Two subgroups of an arbitrary order $p^m$ of $G$
are isomorphic if and only if their types determine the same
partition $(m_1,m_2,...,m_k)$ of $m$ ($0\leq m_1\leq m_2\leq\cdots
\leq m_k$). In this way, we can identify every class of isomorphic
subgroups of $G$ with an element of the direct product
$$C=\xmare{i=1}{k}C_{\alpha_i},$$where $C_{\alpha_i}$ is the chain
$0<1<2<\cdots<\alpha_i$, for all $i=1,2,...,k$. Clearly, given
$[H],[K]\in {\rm Iso}(G)$ ($|H|=p^m$, $|K|=p^{m'}$) that correspond to the partitions
$(m_1,m_2,...,m_k)$ and $(m'_1,m'_2,...,m'_k)$ of $m$ and $m'$,
respectively, we have
$${\rm inf}\{[H],[K]\}=[S] \hspace{1mm}\mbox{ and }\hspace{1mm} {\rm sup}\{[H],[K]\}=[T],$$where
$[S]$ and $[T]$ are determined by the $k$-tuples $({\rm
min}\{m_1,m'_1\},\hspace{1mm} {\rm min}\{m_2,m'_2\},\newline
...,\hspace{1mm}{\rm min}\{m_k,m'_k\})$ and $({\rm
max}\{m_1,m'_1\},\hspace{1mm}{\rm
max}\{m_2,m'_2\},...,\hspace{1mm}{\rm max}\{m_k,m'_k\})$.

Hence ${\rm Iso}(G)$ forms a lattice that can be embedded into a
distributive lattice, namely $C$. This completes the proof.
\hfill\rule{1,5mm}{1,5mm}
\bigskip

Another remarkable class of finite groups for which a similar
conclusion holds is constituted by the so-called {\rm ZM}-groups,
that is the finite groups all of whose Sylow subgroups are cyclic.
It is well-known that two subgroups of such a group $G$ are
conjugate if and only if they have the same order. Moreover,
$C(G)$ is isomorphic to the lattice $L_n$ of all divisors of
$n=|G|$ (see Theorem A of \cite{4}). Thus we infer the following
result.

\bk\n{\bf Proposition 2.2.} {\it Let $G$ be a {\rm ZM}-group of
order $n$. Then the following lattice isomorphisms hold $${\rm
Iso}(G)\cong C(G)\cong L_n.$$In particular, ${\rm Iso}(G)$ is a
distributive lattice.}

\bk\n{\bf Remarks.}
\begin{itemize}
\item[\rm 1.] The conclusion of Proposition 2.2 is valid for the finite cyclic groups,
which are in fact the simplest {\rm ZM}-groups.
\item[\rm 2.] The lattice isomorphism from ${\rm Iso}(G)$ to
$L_n$ in the above proposition is given by $[H]\mapsto |H|$,
for all $[H]\in {\rm Iso}(G)$. For an arbitrary finite
group $G$ of order $n$, this map is only isotone (that is,
$[H]\leq [K]$ implies that $|H|$ divides $|K|$).
\item[\rm 3.] There are finite groups $G$ such that ${\rm Iso}(G)$
is a lattice, but not a distributive or even a modular lattice.
For example, in ${\rm Iso}(S_{2^n})$, $n\geq 4$, the classes
determined by $S_{2^n}$, the three maximal subgroups of $S_{2^n}$
(that are isomorphic to $Q_{2^{n-1}}$, $D_{2^{n-1}}$ and
$\mathbb{Z}_{2^{n-1}}$, respectively) and $\Phi(S_{2^n})$ form a
diamond and so ${\rm Iso}(S_{2^n})$ is not distributive; on the
other hand, we already have seen that ${\rm Iso}(A_4)$ is the
pentagon lattice and therefore it is not modular.
\end{itemize}
\smallskip

Given a finite group $G$, some lattice-theoretical properties can
be transferred from ${\rm Iso}(G)$ to $L(G)$. One of them is the
complementation.

\bk\n{\bf Proposition 2.3.} {\it Let $G$ be a finite group. If
$\hspace{0,5mm}{\rm Iso}(G)$ is a complemented lattice, then
$L(G)$ is also a complemented lattice.}
\bigskip

\n{\bf Proof.} Suppose that ${\rm Iso}(G)$ is a complemented
lattice and denote by $\wedge'$ and $\vee\,'$ its binary
operations. Then for every $[H]\in {\rm Iso}(G)$ there is $[K]\in
{\rm Iso}(G)$ such that $[H]\wedge' [K]=[1]$ and
$[H]\vee\hspace{0,1mm}' [K]=[G]$. Since $H\wedge K$ is contained
both in $H$ and $K$, we infer that $[H\wedge K]\leq [H], [K]$ and
therefore $[H\wedge K]\leq [H]\wedge' [K]$. This implies $H\wedge
K=1$. Similarly, one obtains $H\vee K=G$. Hence $K$ is a
complement of $H$ in $L(G)$.
\hfill\rule{1,5mm}{1,5mm}

\bk\n{\bf Remark.} The converse of Proposition 2.3 is in general
not true. For e\-xam\-ple, $H=\langle y\rangle$ has a complement
in $L(D_8)$, namely $K=\langle x^2, xy\rangle$, but $[K]$ is not a
complement of $[H]$ in ${\rm Iso}(D_8)$ (more precisely, $[H]$ has
no complement in ${\rm Iso}(D_8)$).
\bigskip

Next, we will study the poset of classes of isomorphic subgroups
for finite dihedral groups. Recall that the dihedral group
$D_{2n}$ $(n\ge2)$ is the symmetry group of a regular polygon with
$n$ sides and it has the order $2n$. The most convenient abstract
description of $D_{2n}$ is obtained by using its generators: a
rotation $x$ of order $n$ and a reflection $y$ of order $2$. Under
these notations, we have
$$D_{2n}=\langle x,y\mid x^n=y^2=1,\ yxy=x^{-1}\rangle.$$It is well-known that for every divisor $r$ of $n$, $D_{2n}$ possesses a subgroup isomorphic to $\Z_r$, namely $H^r_0=\langle x^{\frac nr}\rangle$, and $\frac nr$ subgroups isomorphic to $D_{2r}$, namely $H^r_i=\langle x^{\frac nr},x^{i-1}y\rangle,$ $i=1,2,...,\frac nr\hspace{1mm}.$ Then
$$|L(D_{2n})|=\tau(n)+\sigma(n),$$where $\tau(n)$ and $\sigma(n)$ are the number and the sum of all divisors of $n$,
respectively. We easily infer that
$$|{\rm Iso}(D_{2n})|=\left\{\barr{lll}
2\tau(n)-1, \mbox{ for } n \mbox{ odd}\\
&&\\
2\tau(n), \mbox{ for } n \mbox{ even}.\earr\right.$$
\smallskip

We are now able to determine the positive integers $n$ such that
${\rm Iso}(D_{2n})$ forms a lattice.

\bk\n{\bf Proposition 2.4.} {\it The poset $\hspace{0,5mm}{\rm
Iso}(D_{2n})$ is a lattice if and only if either $n$ is odd or
$n=2^k$ for some $k\in\mathbb{N}$.}
\bigskip

\n{\bf Proof.} Suppose first that
$n=p_1^{\alpha_1}p_2^{\alpha_2}\cdots p_k^{\alpha_k}$ with $p_i>2$
prime, $i=1,2,...,k$, or $n=2^k$ for some $k\in\mathbb{N}$. Then,
by using the above description of the subgroups of $D_{2n}$, a
standard induction argument on $k$ easily shows that ${\rm
Iso}(D_{2n})$ forms a lattice.

Conversely, assume that ${\rm Iso}(D_{2n})$ is a lattice, but $n$
is of the form $n=2^{\alpha}\beta$, where $\alpha\geq 1$ and
$\beta \neq 1$ is odd. Then $D_{2n}$ possesses the subgroups
$H=\langle x^{2^{\alpha}},y\rangle\hspace{0,5mm}\cong D_{2\beta}$
and $K=\langle
x^{2^{\alpha-1}}\rangle\hspace{0,5mm}\cong\mathbb{Z}_{2\beta}$.
Clearly, both $H$ and $K$ contain cyclic subgroups of orders 2 and
$\beta$, which proves that $[C_2]\leq [H], [K]$ and
$[C_{\beta}]\leq [H], [K]$ (here $C_2$ and $C_{\beta}$ are
arbitrary cyclic subgroups of $D_{2n}$ of orders 2 and $\beta$,
respectively). It follows immediately that ${\rm inf}\{[H],[K]\}$
does not exist, a contradiction. \hfill\rule{1,5mm}{1,5mm}

\bk\n{\bf Remarks.}
\begin{itemize}
\item[\rm 1.] 12 is the smallest positive integer $n$ for
which there is a finite group $G$ of order $n$ such that ${\rm
Iso}(G)$ is not a lattice.
\item[\rm 2.] Another example of a group with
the above property is the direct pro\-duct
$D_8\times\mathbb{Z}_4$. In this case the classes determined by
the subgroups $H=\langle x^2\rangle\times\hspace{1mm}\mathbb{Z}_4$
and $K=D_8$ does not possess an infimum, too.
\end{itemize}
\smallskip

Unfortunately, we failed in describing exhaustively the class of
finite groups $G$ for which the poset ${\rm Iso}(G)$ is a
(distributive/modular) lattice. We remark that such a group $G$
satisfies the following interesting property: "for every two
distinct prime divisors $p$ and $q$ of $|G|$, either all subgroups
of order $pq$ in $G$ are cyclic or all subgroups of order $pq$ in
$G$ are non-abelian".
\bigskip

We end this section by characterizing all finite groups whose
posets of classes of isomorphic subgroups are chains (i.e.
distributive lattices of a very particular type).

\bk\n{\bf Theorem 2.5.} {\it Let $G$ be a finite group. Then
$\hspace{0,5mm}{\rm Iso}(G)$ is a chain if and only if $G$ is
either a cyclic $p$-group, an elementary abelian $p$-group, a
non-abelian $p$-group of order $p^3$ and exponent $p$ or a
quaternion group of order\, {\rm 8}.}
\bigskip

\n{\bf Proof.} It is clear that for a cyclic $p$-group, an
elementary abelian $p$-group, a non-abelian $p$-group of order
$p^3$ and exponent $p$ or a quaternion group of order 8, the poset
of classes of isomorphic subgroups forms a chain.

Conversely, suppose that ${\rm Iso}(G)$ is a chain. Then $G$ is a
$p$-group. Put $|G|=p^n$ and take a minimal normal subgroup $H$ of
$G$.

If $H$ is the unique subgroup of order $p$\, of $G$, then (4.4) of
\cite{13}, II, shows that $G$ is either cyclic or a generalized
quaternion group $Q_{2^n}$, $n\geq 3$. It is well-known that the
isomorphism classes of the maximal subgroups of $Q_{2^n}$ are
$Q_{2^{n-1}}$ and $\mathbb{Z}_{2^{n-1}}$, and therefore ${\rm
Iso}(Q_{2^n})$ is not a chain for $n\geq 4$. This holds only for
$n=3$, that is for the quaternion group $Q_8$.

If $G$ possesses a minimal subgroup $K$ with $K\neq H$, then $HK$
is elementary abelian of order $p^2$. One obtains that there is no
cyclic subgroup of order $p^2$ in $G$, in other words we have
$${\rm exp}(G)=p.$$

Obviously, if $G$ is abelian, then it is an elementary abelian
$p$-group. In the following we will assume that $G$ is not
abelian. Then $p$ is odd and $G$ contains a non-abelian subgroup
of order $p^3$, say $N$ (more precisely, $N$ is isomorphic with
the group $M(p^3)$, described in (4.13) of \cite{13}, II). Let $A$
be an abelian normal subgroup of maximal order of $G$ and set
$|A|=p^a$. If $a\geq 3$ we infer that $A$ has a subgroup $A_1$ of
order $p^3$. It follows that the classes $[N]$ and $[A_1]$ are not
comparable, a contradiction. In this way, we have $a\leq 2$ and so
$a\in\{1,2\}$. By Corollary 2, \cite{13}, I, page 94, we know that
$2n\leq a(a+1)$, which implies $n\leq 3$. This leads to $n=3$ and
hence $G=N$ is a non-abelian $p$-group of order $p^3$ and exponent
$p$, which completes the proof. \hfill\rule{1,5mm}{1,5mm}
\bigskip

In particular, Theorem 2.5 shows that there are only two finite
non-abelian groups $G$ with ${\rm Iso}(G)$ fully ordered, and each
of them is of order $p^3$ for some prime $p$.

\section{Finite groups with the same poset\\ of classes of isomorphic subgroups}

In this section we study when for two finite groups $G_1$ and
$G_2$ the poset/lattice isomorphism ${\rm Iso}(G_1)\cong {\rm
Iso}(G_2)$ holds. Obviously, a sufficient condition to have this
isomorphism is $G_1\cong G_2$, but it is not necessary, as show
the examples in Section 2. We also remark that the weaker
condition $L(G_1)\cong L(G_2)$ does not imply that ${\rm
Iso}(G_1)\cong {\rm Iso}(G_2)$ (for example, take $G_1=S_3$ and
$G_2=\mathbb{Z}_3\times\mathbb{Z}_3$) and the same thing can be
said about the converse implication.
\bigskip

We start with the following easy but important lemma.

\bk\n{\bf Lemma 3.1.} {\it Let $G_1$ be a finite $p$-group of
order $p^n$. If $G_2$ is a finite group such that ${\rm
Iso}(G_1)\cong {\rm Iso}(G_2)$, then $G_2$ is a $q$-group of order
$q^n$.}
\bigskip

\n{\bf Proof.} The condition $|G_1|=p^n$ implies that ${\rm
Iso}(G_1)$ possesses a unique non-trivial element, namely the
class determined by the subgroups of order $p$. Then ${\rm
Iso}(G_2)$ satisfies a similar property, and so $G_2$ is a
$q$-group for some prime $q$. One the other hand, we easily infer
that all maximal chains of ${\rm Iso}(G_1)$ are of length $n$.
Since a poset isomorphism preserves the length of such a chain,
one obtains that $|G_2|=q^n$, as desired.
\hfill\rule{1,5mm}{1,5mm}
\bigskip

The above lemma can be extended to finite groups of arbitrary
orders in the following manner.

\bk\n{\bf Theorem 3.2.} {\it Let $G_1$ and $G_2$ be two finite
groups such that ${\rm Iso}(G_1)\cong {\rm Iso}(G_2)$. If\,
$|G_1|=p_1^{\alpha_1}p_2^{\alpha_2}\cdots p_k^{\alpha_k}$, where
$p_i$, $i=1,2,...,k$, are distinct primes, then we have
$|G_2|=q_1^{\alpha_1}q_2^{\alpha_2}\cdots q_k^{\alpha_k}$ for some
distinct primes $q_1,q_2,...,q_k$.}
\bigskip

\n{\bf Proof.} Let $f:{\rm Iso}(G_1)\longrightarrow {\rm
Iso}(G_2)$ be a poset isomorphism, $i\in\{1,2,...,k\}$ and $S_i$
be a Sylow $p_i$-subgroup of $G_1$. If $f([S_i])=[S_i']$, then, by
Lemma 3.1, we have $|S_i'|=q_i^{\alpha_i}$ for some prime
$q_i$. Moreover, it is easy to see that $S_i'$ is a Sylow subgroup
of $G_2$. Hence $|G_2|$ is of the form
$q_1^{\alpha_1}q_2^{\alpha_2}\cdots q_k^{\alpha_k}$ with
$q_1,q_2,...,q_k$ distinct primes, completing the proof.
\hfill\rule{1,5mm}{1,5mm}
\bigskip

Proposition 2.1 and Theorem 3.2 lead the following immediate
characterization of the poset isomorphism ${\rm Iso}(G_1)\cong
{\rm Iso}(G_2)$ for two finite abelian
groups $G_1$ and $G_2$.

\bk\n{\bf Corollary 3.3.} {\it Let $G_1\hspace{-1mm}$ and $G_2$ be
two finite abelian groups of orders
$p_1^{\alpha_1}\hspace{-1mm}p_2^{\alpha_2}\hspace{-1mm}\cdots\hspace{-1mm}
p_k^{\alpha_k}$ and $q_1^{\beta_1}q_2^{\beta_2}\cdots
q_r^{\beta_r}$, respectively. Then $\hspace{0,5mm}{\rm Iso}(G_1)\cong{\rm Iso}(G_2)$ if and only if $k=r$\newline and there is a permutation $\sigma$ of $\{1,2,...,k\}$ such that $\beta_i=\alpha_{\sigma(i)}$ and $\hspace{0,5mm}{\rm Iso}(S_{q_i}')\cong {\rm Iso}(S_{p_{\sigma(i)}})$, where $S_{q_i}'$ is the Sylow $q_i$-subgroup of $G_2$ and $S_{p_{\sigma(i)}}$ is the Sylow $p_{\sigma(i)}$-subgroup of $G_1$, for all $i=1,2,...,k$.}

\bk\n{\bf Example.} By using Corollary 3.3, we easily infer that $${\rm Iso}(\mathbb{Z}_2\times\mathbb{Z}_6\times\mathbb{Z}_{18})\cong {\rm
Iso}(\mathbb{Z}_7\times\mathbb{Z}_{6125}).$$

Lemma 3.1 shows that the class of finite $p$-groups is preserved
by isomorphisms between their posets of classes of isomorphic
subgroups. Other important classes of finite groups satisfying the
same property are solvable groups, {\rm CLT}-groups and
supersolvable groups, respectively. This is due to the fact that
if $f:{\rm Iso}(G_1)\longrightarrow {\rm Iso}(G_2)$ is a poset
isomorphism and $G_1$ is of one of the above three types, then the
strong connection between the orders of $G_1$ and $G_2$ assures
the existence of Hall subgroups, of subgroups of any orders or the
validity of the Jordan-Dedekind chain condition for $G_2$.

\bk\n{\bf Corollary 3.4.} {\it The classes of finite solvable
groups, {\rm CLT}-groups and supersolvable groups are preserved by
isomorphisms between their posets of classes of isomorphic
subgroups.}
\bigskip

We end this section by some results related to the uniqueness of a
finite group with a given poset of classes of isomorphic
subgroups.

\bk\n{\bf Theorem 3.5.} {\it Let $n\geq 2$ be an integer which is
not square-free. Then there are at least two non-isomorphic groups $G_1$
and $G_2$ of order $n$ such that ${\rm Iso}(G_1)\cong {\rm
Iso}(G_2)$.}
\bigskip

\n{\bf Proof.} Let $n=p_1^{\alpha_1}p_2^{\alpha_2}\cdots
p_k^{\alpha_k}$ be the decomposition of $n$ as a product of prime
factors. We will proceed by induction on $k$.

Suppose first that $k=1$, that is $n=p_1^{\alpha_1}$ with
$\alpha_1\geq 2$. For $\alpha_1=2$ we take
$G_1=\mathbb{Z}_{p_1^2}$ and
$G_2=\mathbb{Z}_{p_1}\times\mathbb{Z}_{p_1}$. For $\alpha_1\geq 3$
we take $G_1=M(p_1^{\alpha_1})$ (see Theorem 4.1 of \cite{13}, II)
and $G_2=\mathbb{Z}_{p_1^{\alpha_1-1}}\times\mathbb{Z}_{p_1}$ if
$p_1\neq 2$, respectively $G_1=D_{p_1^{\alpha_1}}$ and
$G_2=\mathbb{Z}_{p_1^{\alpha_1-1}}\times\mathbb{Z}_{p_1}$ if
$p_1=2$.

Assume now that $k\geq 2$. By the inductive hypothesis, we can
choose two non-isomorphic groups $H_1$ and $H_2$ of the order
$p_1^{\alpha_1}p_2^{\alpha_2}\cdots p_{k-1}^{\alpha_{k-1}}$ such
that ${\rm Iso}(H_1)\cong {\rm Iso}(H_2)$. Then it is a simple
exercise to see that the groups
$G_1=H_1\times\mathbb{Z}_{p_k^{\alpha_k}}$ and
$G_2=H_2\times\mathbb{Z}_{p_k^{\alpha_k}}$ satisfy the desired
conditions, completing the proof.
\hfill\rule{1,5mm}{1,5mm}
\bigskip

Inspired by Theorem 3.5, we came up with the following conjecture,
which we have verified for several finite groups of small orders.

\bk\n{\bf Conjecture.} {\it For every non-trivial finite group
$G_1$ whose order is not square-free there exists a finite group $G_2$
such that $|G_1|=|G_2|$, ${\rm Iso}(G_1)\cong {\rm Iso}(G_2)$ and
$G_1\ncong G_2$.}

\bk\n{\bf Remark.} The above conjecture says nothing else than the
implication
$$|G_1|=|G_2|\mbox{ and } {\rm Iso}(G_1)\cong {\rm
Iso}(G_2)\hspace{1mm}\Longrightarrow\hspace{1mm} G_1\cong
G_2$$fails in all cases, except when $|G_1|=|G_2|$ is a
square-free number.

\section{Conclusions and further research}

All our previous results show that the poset consisting of all
classes of isomorphic subgroups of a (finite) group can constitute
a significant aspect of (finite) group theory. Clearly, the study
started in this paper can successfully be extended to other
classes of groups. It can be also generalized by studying the
posets of isomorphic substructures of other algebraic structures
(rings, modules, algebras, ... and so on). This will surely be the
subject of some further research.
\bigskip

Finally, we mention several open problems concerning this topic.
\bigskip

\noindent{\bf Problem 4.1.} Determine the finite groups $G$ for
which the poset Iso($G$) is a lattice and study the properties of
this lattice.
\bigskip

\noindent{\bf Problem 4.2.} What can be said about two
\textit{arbitrary} finite groups $G_1$ and $G_2$ satisfying ${\rm
Iso}(G_1)\cong {\rm Iso}(G_2)$?
\bigskip

\noindent{\bf Problem 4.3.} Given two finite groups $G_1$ and
$G_2$, study the isomorphisms between the posets/lattices ${\rm
Iso}(G_1)$ and ${\rm Iso}(G_2)$ induced by the isomorphisms or by
the $L$-isomorphisms between $G_1$ and $G_2$.
\bigskip

\noindent{\bf Problem 4.4.} The most natural generalization of the
poset Iso($G$) asso\-cia\-ted to a finite group $G$ is obtained by
considering $L$-isomorphisms instead of group isomorphisms in its
definition: $${\rm Iso}'(G)=\{[H]' \mid H\in L(G)\}, \mbox{ where
} [H]'=\{K\in L(G) \mid L(K)\cong L(H)\}.$$Investigate the above
new poset ${\rm Iso}'(G)$ with respect to the same ordering
relation as for Iso($G$).
\bigskip

\noindent{\bf Problem 4.5.} Given a finite group $G$, study the
posets of classes of subgroups with respect to other equivalence
relations on $L(G)$ (or on other important subposets of $L(G)$).
For example:
\begin{itemize}
\item[1.] $H\sim_1 K$ if and only if $|H|=|K|$\,;
\item[2.] $H\sim_2 K$ if and only if there is $f\in {\rm Aut}(G)$ such that $f(H)=K$;
\item[3.] $H\sim_3 K$ if and only if $\pi_e(H)=\pi_e(K)$ (that is, $H$ and $K$ have the same set of element orders).
\end{itemize}
\bigskip

\noindent{\bf Problem 4.6.} The concept of \textit{solitary
quotient} of a finite group has been defined in \cite{15} as the
dual concept of \textit{solitary subgroup}. Following the same
technique, we can construct a "dual" for the set Iso($G$), namely
$${\rm QIso}(G)=\{[H] \mid H\in N(G)\}, \mbox{ where }
[H]=\{K\in N(G) \mid G/K\cong G/H\}.$$Endow this set with a
suitable ordering relation and study some similar problems.
\bigskip

\bigskip\noindent{\bf Acknowledgements.} The author is grateful to the reviewer for
its remarks which improve the previous version of the paper.

\vspace*{5ex}\small

\hfill
\begin{minipage}[t]{5cm}
Marius T\u arn\u auceanu \\
Faculty of  Mathematics \\
``Al.I. Cuza'' University \\
Ia\c si, Romania \\
e-mail: {\tt tarnauc@uaic.ro}
\end{minipage}

\end{document}